\documentclass[11pt]{article}

\usepackage{amsfonts, mathrsfs}
\usepackage{amsmath}
\usepackage{latexsym}

\newtheorem{theorem}{Theorem}
\newtheorem{remark}{Remark}

\begin{document}

\begin{center}
{\LARGE The exponentiated exponential Poisson distribution revisited}\\[1ex]
Tibor K. Pog\'any \\[1ex]
{\it Universitas Budensis, John von Neumann Faculty of Informactics, 1034 Budapest, Hungary}\\
{\it University of Rijeka, Faculty of Maritime Studies, 51000 Rijeka, Croatia}
\footnote{e--mail: tkpogany@gmail.com, poganj@pfri.hr}\\
\end{center}

\begin{abstract}
The main aim of this article is to characterize and investigate the three parameter exponentiated exponential Poisson probability distribution ${\rm EEP}(\alpha, \beta, 
\lambda)$ by giving explicit closed form expressions for its characteristic function $\phi_\xi(t)$ and moment generating function $M_\xi(t)$, and finally, to show that the 
existing series and integral form expressions for positive integer order moments 
$\mathbb E\xi^\nu, \nu \in \mathbb N$ are in fact valid for all $\nu>1-\alpha, \alpha>0$. 
\medskip

\medskip

\end{abstract}

\section{Introduction}

When we start from  a composite series connected system, which consists form several blocks built from $\alpha$ independent identically distributed units in paralellel connected, then the time to failure $\xi$ becomes that, of parallel blocks 
in series. According to Risti\'c and Nadarajah \cite{RisNad1} the starting units lifetime distribution being exponential $\mathscr E(\beta)$ and the number $N$ of parallel--blocks is distributed according to {\em zero--truncated} Poisson $\mathscr P_0(\lambda)$ law \cite{Kus}, that is
   \[ \mathbb P\{N=n\} = \frac{\lambda^n}{n!\,(1-\exp(-\lambda))}, \qquad n \in \mathbb N\, ,\]
the unconditonal {\it cumulative distribution function} (CDF) of the time to failure of the first out of the $N$ functioning blocks becomes an {\em exponentiated exponential Poisson} distributed rv $\xi$, having CDF \cite[Eq. (2)]{RisNad1}
   \begin{equation} \label{T1}
      F_\xi(x) = \frac1{1-{\rm e}^{-\lambda}}\,
             \left(1-{\rm e}^{-\lambda(1-{\rm e}^{-\beta x})^\alpha}\right)\,
             \chi_{\mathbb R_+}(x),
   \end{equation}
where $\theta = (\alpha, \beta, \lambda)>0$. Here $\chi_S(x)$ stands for the characteristic function of the set $S$, that is $=1$ if $x \in S$ and $=0$ else. Thus, this type rv considered on a given standard probability space $(\Omega, \mathscr F, \mathbb P)$ we write throughout 
$\xi \sim {\rm EEP}(\theta)$. The related {\it probability density function} (PDF) \cite[Eq. (2)]{RisNad1}
   \begin{equation} \label{T2}
      f_\xi(x) = \frac{\alpha\beta\lambda}{1-{\rm e}^{-\lambda}}\,
             {\rm e}^{-\beta x}\,\left(1-{\rm e}^{-\beta x}\right)^{\alpha-1}\,
             {\rm e}^{-\lambda(1-{\rm e}^{-\beta x})^\alpha}\, \chi_{\mathbb R_+}(x),
   \end{equation}
with $\beta$ a scaling, and both $\alpha, \lambda$ shape parameters. 

Risti\'c and Nadarajah \cite{RisNad1} give four main reasons for introducing EEP($\theta$), 
such as: {\bf (i)} the previously explaned failure rate's mathematical model; {\bf (ii)} EEP($\theta$) distribution exhibits monotone increasing, monotone decreasing and upside--down bathtub hazard rate's behaviour; {\bf (iii)} being $F(x) \sim \mathscr O(x^\alpha), x \to 0$ and $F(x) \sim \mathscr O({\rm e}^{-\beta x}), x \to \infty$ (compare \eqref{T1} for both cases), EEP($\theta$) has the lower tail  behaviour like the gamma, Weibull and allied distributions, and the upper tail behaviour indentical to the exponential distribution; {\bf (iv)} the EEP($\theta$)  is superior in fitting rel data with respect to a wide selection (in \cite{RisNad1} have been tested 15) of the most frequent two--, or three--parameter distributions. Finally, the reader's attention is drawn for further informations to the ancestor article \cite{RisNad1} and the related references therein. 

Our motivation here is threefold: {\bf 1.} to give a closed form expression for the {\it characteristic function} (CHF) in terms of certain special functions; {\bf 2.} to establish a closed expression for the {\it moment generating functon} (MGF) and {\bf 3.} to give a proof that the existing series representation of the positive integer order moments by Risti\'c and Nadarajah \cite[p. 6, Eq. (7)]{RisNad1} holds true for a real order as well. 

It is also worth to mention that as an immediate step, a closed form expression was obtained  for the general order moment in the case of {\em exponentiated exponential} ${\rm EE}(\alpha, \beta)$  distribution, see section III.

In the derivation procedures we use generalized hypergeometric type special functions like Fox--Wright ${}_1\Psi_1$ and the Goyal--Laddha generalized Hurwitz--Lerch Zeta $\Phi^*_\mu$.

\section{The characteristic function for ${\rm EEP}(\alpha, \beta, \lambda)$} 

The Fox--Wright generalization ${}_p\Psi_q$ of the generalized hypergeometric function ${}_pF_q$, with $p$ numerator and $q$ denominator parameters reads \cite{SGG}
   \[ {}_p\Psi_q \Big[ \begin{array}{ccc} (\alpha_1, A_1), \cdots, (\alpha_p, A_p)\\ 
         (\beta_1, B_1), \cdots, (\beta_q, B_q) \end{array} \Big|\, x\Big] 
       = \sum_{n=0}^\infty \frac{\prod\limits_{k=1}^p 
         \Gamma \big(\alpha_k + A_k\, n\big)}{\prod\limits_{k=1}^q 
         \Gamma \big(\beta_k+B_k\, n\big)}\,\frac{x^n}{n!}\, , \]
where $\alpha, \beta \in \mathbb C$, $A, B>0$; and the series converges for suitably bounded values of $|x|$ when $\Delta = 1+ \sum_{k=1}^q B_k - \sum_{j=1}^p A_j>0$, and an empty product  is by convention used to be 1. In t{ur}n, ${}_1\Psi_1$ we call {\it confluent}.

The characteristic function of a rv $\xi$ having PDF $f_\xi(x)$ and CDF $F_\xi(x)$ is  actually the Fourier transform of the PDF (or the Fourier--Stieltjes transform of the CDF):
   \[ \phi_\xi(t) = \mathbb E\big({\rm e}^{{\rm i}t\xi}\big) 
                  = \int_{\mathbb R} {\rm e}^{{\rm i}tx}\,f_\xi(x)\, {\rm d}x
                  = \int_{\mathbb R} {\rm e}^{{\rm i}tx}\,{\rm d}F_\xi(x);\]
the related moment generating function we define as
   \[ M_\xi(t) = \phi_\xi({\rm i}t)\, .\]
Now, we formulate closed form expression results for both CHF and MGF for the exponentiated 
exponential Poisson distribution. 

\begin{theorem} {\it Let $\xi \sim {\rm EEP}(\theta)$, $\theta = (\alpha, \beta, \lambda)>0$. For all $t \in \mathbb R$ we have
   \begin{equation} \label{U2}
      \phi_\xi(t) = \frac{\lambda \Gamma\left(1-\frac{{\rm i}t}\beta \right)}
         {1-{\rm e}^{-\lambda}}\, \,{}_1\Psi_1\Big[ \begin{array}{ccc} (\alpha, \alpha)\\
         (1+\alpha-\frac{{\rm i}t}\beta, \alpha)\end{array} \Big|\, -\lambda \,\Big]\, .
   \end{equation}
Moreover, for all $t>-\beta$}
   \begin{equation} \label{U3}
      M_\xi(t) = \frac{\lambda \Gamma\left(1+\frac{t}\beta \right)}
           {1-{\rm e}^{-\lambda}}\, {}_1\Psi_1\Big[ \begin{array}{ccc} (\alpha, \alpha)\\ 
           (1+\alpha+\frac{t}\beta,\alpha) \end{array} \Big|\, -\lambda \,\Big]\, .
   \end{equation}
\end{theorem}

\noindent {\it Proof}: Having in mind \eqref{T2} we conclude
   \[\phi_\xi(t) = C \int_0^\infty {\rm e}^{{\rm i}tx - \beta x}
                   \left(1-{\rm e}^{-\beta x}\right)^{\alpha-1}
                   {\rm e}^{-\lambda(1-{\rm e}^{-\beta x})^\alpha} \,{\rm d}x;\]
this, by substitution $1-\exp(-\beta x) = u$, yields the {\em mutatis mutandis} result
   \begin{equation} \label{U4}
      \phi_\xi(t) = \frac{\alpha\lambda}{1-{\rm e}^{-\lambda}}\,\int_0^1 
                    (1-u)^{-\frac{{\rm i}t}\beta}\, u^{\alpha-1}\,
                    {\rm e}^{-\lambda u^\alpha}\, {\rm d}u\, .
   \end{equation}  
The Maclaurin expansion of the exponential term and interchanging the order of summation and 
integration, lead us to  
   \[ \phi_\xi(t)= \frac{\alpha\lambda\,\Gamma(1-\frac{{\rm i}t}\beta)}
                   {1-{\rm e}^{-\lambda}}\, \sum_{n \geq 0} 
                   \frac{\Gamma(\alpha + \alpha n)}{\Gamma(1-\frac{{\rm i}t}\beta + 
                   \alpha + \alpha n)}\, \frac{(-\lambda)^n}{n!}\, ,\]
which is actually \eqref{U2}. Setting $t=0$ in \eqref{U2}, we have to have $\phi_\xi(0) = 1$. Indeed,
   \[ \phi_\xi(t) = \frac{\alpha\lambda}{1-{\rm e}^{-\lambda}}\, 
                     \sum_{n \geq 0}\frac{\Gamma(\alpha + \alpha n)}
                     {\Gamma(1+\alpha+\alpha n)}\, \frac{(-\lambda)^n}{n!} 
                  = \frac{\lambda}{1-{\rm e}^{-\lambda}}\, 
                     \sum_{n \geq 0} \frac{(-\lambda)^n}{n!\,(n+1)} = 1\,.\]
The rest is clear. 
\hfill $\Box$

\begin{remark} There is an analytic continuation of $M_\xi(t)$ from $(-\beta, \infty)$ to the 
whole $\mathbb C \setminus \bigcup_{n \geq 1}(\lfloor{-\beta\rfloor} - n)$, where $\lfloor{x}\rfloor$ denotes the integer part of $x$. Namely, all poles of $M_\xi$ are simple and occur at $-\beta m, m \in \mathbb N$. 
\end{remark}

\begin{remark} An integral representation, similar to \eqref{U4}, appears in \cite[p. 5]{RisNad1}, where erroneously stands $\theta$ instead of the correct variable $t$. However, the authors didn't reach the formula \eqref{U3} for the MGF $M_\xi(t)$. 
\end{remark}

\section{The moment properties}

Consider the rv $\eta$ on a fixed probability space $(\Omega, \mathscr F, \mathbb P)$ having CDF
   \[ F_\eta(x) = \left(1-{\rm e}^{-\beta x}\right)^\alpha \chi_{\mathbb R_+}(x), 
                  \qquad \alpha, \beta>0\, .\]
Then $\eta$ behaves according the {\em exponentiated exponential distribution} ${\rm EE}(\alpha, \beta)$ which, due among others Gupta and Kundu \cite{GupKun1, GupKun2}, attracts notable attention. We point out that the Maclaurin series expansion of the CDF 
   \begin{equation} \label{R0}
      F_\xi(x) = \frac1{1-{\rm e}^{-\lambda}} \sum_{m \geq 1} 
                 \frac{(-1)^{m+1}}{m!}\,[\lambda\,F_\eta(x)]^m
   \end{equation}
clearly shows that ${\rm EEP}(\theta)$ is a mixture of ${\rm EE}(\alpha, \beta)$ in the sense that the associated parameters coincide.

According to \cite[p. 488, Eq. (1.1)]{SSPS} we introduce the Hurwitz--Lerch Zeta (HLZ) function 
   \[ \Phi(z, s, a) = \sum_{n \geq 0} \frac{z^n}{(n+a)^s}\, ,\]
where $a \in \mathbb C \setminus \mathbb Z_0^-$, $s \in \mathbb C$ when $|z|<1$; $\Re s >1$ 
for $|z|=1$. Further, the Goyal--Laddha generalized HLZ function \cite[p. 100, Eq. (1.5)]{GoyLad}
   \begin{equation} \label{R1}
      \Phi_\mu^*(z, s, a) = \sum_{n \geq 0} \frac{(\mu)_n}{n!}\,\frac{z^n}{(n+a)^s}\, ,
   \end{equation}
where $\mu \in \mathbb C$; $a \in \mathbb C \setminus \mathbb Z_0^-$, $s \in \mathbb C$ when $|z|<1$; $\Re (s-\mu)>1$ for $|z|=1$. We also list the integral form \cite[p. 100, Eq. (1.6)]{GoyLad}, \cite[p. 495, Eq. (2.10)]{SSPS}
   \begin{equation} \label{R2}
      \Phi_\mu^*(z, s, a) = \frac1{\Gamma(s)} \int_0^\infty \frac{t^{s-1}\,{\rm e}^{-st}}
                            {(1-z{\rm e}^{-t})^\mu}\, {\rm d}t\, .
   \end{equation}
Gupta and Kundu \cite{GupKun1, GupKun2} showed that
   \[ \mathbb E\eta^n = \frac{\alpha\,n!}{\beta^n}
               \sum_{k \geq 0} \frac{(-1)^k}{(k+1)^{n+1}} \binom{\alpha -1}k, 
               \quad k \in \mathbb N\, ,\]
where $n \in \mathbb N$, and the generalized binomial coefficient 
$\tbinom wk = (-1)^k\,(-w)_k\,(k!)^{-1},\, w \in \mathbb C$. The related closed form representation results follow. 

\begin{theorem} {\it Let $\eta$ be a ${\rm EE}(\alpha, \beta)$ rv, $\alpha, \beta>0$. Then for all $\nu>1-\alpha$ we have}
   \begin{equation} \label{R5}
      \mathbb E\eta^\nu = \frac{\alpha\,\Gamma(\nu+1)}{\beta^\nu}
               \Phi_{1-\alpha}^*(1, \nu+1, 1)\, .
   \end{equation}
\end{theorem}

\noindent {\it Proof}: Obviously, being $\beta>0$, expressing the binomial in the integrand as a power series, we conclude
   \begin{align*}
      \mathbb E\eta^\nu &= \int_0^\infty x^\nu\,F_\eta'(x)\,{\rm d}x 
             = \alpha\beta \int_0^\infty x^\nu\,{\rm e}^{-\beta x}
               \left(1-{\rm e}^{-\beta x}\right)^{\alpha-1}\,{\rm d}x\\
            &= \alpha\beta \sum_{k \geq 0} (-1)^k \binom{\alpha-1}k \int_0^\infty 
               x^\nu\,{\rm e}^{-\beta (k+1)x}\,{\rm d}x\, .
   \end{align*}
Applying the known Gamma function formula
   \[ \Gamma(a) = A^a \int_0^\infty x^{a-1}{\rm e}^{-Ax}\,{\rm d}x, 
                \quad \min\{\Re(a), \Re(A)\}>0\, ,\]
we get
   \[ \mathbb E\eta^\nu = \frac{\alpha\,\Gamma(\nu+1)}{\beta^{\nu}}\, \sum_{k \geq 0}
            \frac{(-1)^k}{(k+1)^{\nu+1}}\, \frac{(\alpha-1)\cdots(\alpha-k)}{k!} \, .\]         
By the Pochhammer symbol property $a(a-1)\cdots(a-k+1) = (-1)^k (-a)_k, k \in \mathbb N_0$ 
we infer
   \[ \mathbb E\eta^\nu = \frac{\alpha\,\Gamma(\nu+1)}{\beta^{\nu}}\, \sum_{k \geq 0}
            \frac{(1-\alpha)_k}{(k+1)^{\nu+1}} \frac1{k!}\, ,\]
which, in conjunction with the definition \eqref{R1} of the Goyal--Laddha HLZ function gives \eqref{R5}. Since $\Delta = \Re(\nu+1-(1-\alpha)) = \nu+\alpha>1$, the proof is completed.
\hfill $\Box$

Risti\'c and Nadarajah \cite[p. 6, Eq. (10)]{RisNad2} reported on the following double series  representation result:
   \begin{equation} \label{R4}
      \mathbb E\xi^n = \frac{\alpha \lambda \,n!}{\beta^n(1-{\rm e}^{-\lambda})}
               \sum_{m, k \geq 0} \frac{(-1)^{m+k}\lambda^m}{m!\,(k+1)^{n+1}}\,
               \binom{\alpha(m+1)-1}k\,. 
   \end{equation}
Now, we show that \eqref{R4} can be extended to general order 
moments case, reducing the double sum expression into a simple one, establishing the main conclusion on the widest possible parameter range. 

\begin{theorem} {\it Let $\xi \sim {\rm EEP}(\theta)$, $\theta = (\alpha, \beta, \lambda)>0$. For all $\nu>1-\alpha$ there holds true}
   \begin{equation} \label{R6}
      \mathbb E\xi^\nu = \frac{\alpha \lambda \,\Gamma(\nu+1)}{\beta^\nu (1-{\rm e}^{-\lambda})}
                         \sum_{m \geq 0} \Phi^*_{1-\alpha(m+1)}(1, \nu+1, 1)\,
                         \frac{(-\lambda)^m}{m!}. 
   \end{equation}
\end{theorem}

\noindent {\it Proof}: Starting with
   \[ \mathbb E\xi^\nu = \int_0^\infty x^\nu\,f_\xi(x)\, {\rm d}x 
                        = \frac{\alpha\lambda}{\beta^{\nu}\big(1-{\rm e}^{-\lambda}\big)}\,
                   \int_0^\infty x^\nu {\rm e}^{-x}
                   \left(1 - {\rm e}^{-x}\right)^{\alpha-1}
                   {\rm e}^{-\lambda(1-{\rm e}^{-x})^\alpha}{\rm d}x,\]
then expanding $\exp\{-\lambda(1-{\rm e}^{-x})^\alpha\}$ into a Maclaurin series and 
recalling the Goyal--Laddha HLZ function's integral expression \eqref{R2} we deduce
   \begin{align*}
      \mathbb E\xi^\nu &= \frac{\alpha\lambda}{\beta^{\nu}\big(1-{\rm e}^{-\lambda}\big)}\,
                   \sum_{m \geq 0} \frac{(-\lambda)^m}{m!}
                   \int_0^\infty \frac{x^\nu {\rm e}^{-x}}
                   {(1 - {\rm e}^{-x})^{1-\alpha(m+1)}}\,{\rm d}x\\
                &= \frac{\alpha\lambda\Gamma(\nu+1)}{\beta^{\nu}
                   \big(1-{\rm e}^{-\lambda}\big)}\, 
                   \sum_{m \geq 0} \frac{(-\lambda)^m}{m!}
                   \Phi_{1-\alpha(m+1)}^*(1, \nu+1, 1)\, .
   \end{align*}
Since the convergence issues $\Delta = \Re(\nu+1-(1-\alpha(m+1))= \nu+\alpha(m+1) \geq 
\nu + \alpha>1,\,m \in \mathbb N$ are satisfied by convention \eqref{R6} is proved. 
\hfill $\Box$

\begin{remark}
The same conclusion can be made by applying \eqref{R0}, that is
   \[ \mathbb E\xi^\nu = \frac1{1-{\rm e}^{-\lambda}} \sum_{m \geq 1} 
               \frac{(-1)^{m-1}\lambda^m}{m!}\,\int_0^\infty x^\nu\,d[F_\eta(x)]^m\, .\]
Now, subsequent suitable transformations lead to \eqref{R6}. 
\end{remark} 

\section{Concluding notes}

\noindent {$\mathsf A.$} Closed form expressions are established, for the first time, for the characteristic function $\phi_\xi(t) = \mathbb E{\rm e}^{{\rm i}t\xi}$ and the moment generating function $M_\xi(t) = \phi_\xi({\rm i}t)$ in the case of exponentiated exponential Poisson distribution $\xi \sim {\rm EEP}(\theta)$, $\theta>0$ generalizing the results by Risti\'c and Nadarajah \cite{RisNad1, RisNad2}. The findings are presented in terms of the confluent Fox--Wright ${}_1\Psi_1$ function as Theorem 1. \medskip

\noindent {$\mathsf B.$} Inspecting \cite{RisNad1} and the references therein, we see that 
the moment problem for a rv coming from ${\rm EEP}(\theta)$ distribution family was 
solved {\em via} double series, and separately in a definite integral form, the both {\em exclusively} for the moments of positive integer order. 

Here, firstly we are obtain a closed form representation for the real order moment in the case of exponentiated exponential distribution ${\rm EE}(\alpha, \beta)$ using Goyal--Laddha  generalized HLZ function. This novel result is exposed in Theorem 2. In the sequel, using similar technique, the real order moments problem has been solved for ${\rm EEP}(\theta)$ distribution, in the form of a simple weighted sum of Goyal--Laddha generalized HLZ functions (Theorem 3). 

It is worth to mention the existence results of negative real order moments in both cases, which are by--products of Theorems 2 and 3. \medskip

\noindent {$\mathsf C.$} There exists the Riemann--Liouville fractional derivation formula 
\cite[p. 490, Eq. (1.17)]{SSPS}
   \[ \Phi_\mu^*(z, s, a) = \frac1{\Gamma(\mu)}\, \mathcal D_z^{\mu-1}\left\{ z^{\mu-1}
                            \Phi(z,s,a)\right\}, \quad \Re(\mu)>0\, ;\]
this exhibits the useful fact that the Goyal--Laddha generalized HLZ function $\Phi_\mu^*(z, s, a)$ is essentially a consequence of the classical Hurwitz--Lerch Zeta function $\Phi(z, s, a)$. So, the results presented in Theorems 2 and 3 can be treated in this manner as well. 
\medskip

\noindent {$\mathsf D.$} Finally, the argumentation by Risti\'c and Nadarajah \cite[Introduction]{RisNad1} and their real data fitting comparation shows that ${\rm EEP}(\theta)$ is an excellent candidate for lifetime distribution setting in reliability problem solving studies. Hence, immediately arises the question of the reliability equivalence analysis \cite{PTT, PTV} of composite series and parallel systems  having independent identically distributed units possessing ${\rm EEP}(\theta)$ lifetime distribution. This could be realized by the so--called reduction method developed by R{\aa}de \cite{R1, R2}; alternatively Sarhan \cite{S1} introduces the reliability equivalence factor by which the failure rates of some of the system's components should be reduced in order to reach equality of the reliability of another better system. R{\aa}de and Sarhan considered components with exponential lifetime distribution, while the gamma--lifetime was treated in \cite{XZ}. 
Comparing the reduction method for ${\rm EEP}(\theta)$ lifetime distribution to the associated hot--, and cold--duplication would be of considerable interest as well. \medskip

However, we leave these questions unanswered here, it belong to some future studies.


\begin{thebibliography}{99}

\bibitem{RisNad1}
Risti\'c MM, Nadarajah S. A new lifetime distribution. {\it Journal of Statistical Computation and Simulation} 2012. [DOI:10.1080/00949655.2012.697163]

\bibitem{Kus}
Ku\c s C. A new lifetime distribution. {\it Comput. Statist. Data Anal.} 2007; 51(9): 4497--4509. 

\bibitem{RisNad2}
Risti\'c MM, Nadarajah S. A new lifetime distribution. Research Report No. 21. 2010; Probability and Statistics Group, School of Mathematics, The University of Manchester.

\bibitem{SGG} 
Srivastava HM, Gupta KC, Goyal SP. The $H$--Functions of One and Two 
Variables with Applications. New Delhi : South Asian Publishers; 1982.

\bibitem{GupKun1}
Gupta RD, Kundu D. Generalized exponential distributions. {\it Aust. N. Z. J. Stat.} 
1999; 41(2): 173--188. 
 
\bibitem{GupKun2}
Gupta RD, Kundu D. Exponentiated exponential family: an alternative to gamma and Weibull distributions. {\it Biom. J.} 2001; 43(1): 117--130.

\bibitem{SSPS}
Srivastava HM, Saxena RK, Pog\'any TK, Saxena R. Integral and computational  representations of the extended Hurwitz–Lerch zeta function, {\it Integral Transforms Spec.  Funct.} 2011; 22(7):  487–-506.

\bibitem{GoyLad}
Goyal SP, Laddha RK. On the generalized Zeta function and the generalized 
Lambert function. {\it Ga\d nita Sandesh} 1997; 11: 99–-108. 

\bibitem{PTT}
Pog\'any TK, Tomas V, Tudor M. Hot duplication versus survivor equivalence in Gamma-Weibull  distribution. {\it J. Stat. Appl. Prob.} 2013; 2(1): 1--10. 

\bibitem{PTV}
Pog\'any TK, Tudor M, Val\v ci\'c S. Cold duplication and survival equivalence in the case of gamma -- Weibull distributed composite systems [submitted manuscript]. 2014. 

\bibitem{R1} 
R{\aa}de L. Reliability equivalence. {\it Microelectronics and Reliability}. 1998; 33: 323--325.
 
\bibitem{R2}
R{\aa}de L. Reliability survival equivalence. {\it Microelectronics and Reliability.} 
1993; 33: 881--894. 

\bibitem{S1} 
Sarhan A. Reliability equivalence with a basic series/parallel system. {\it App. Math. Comput.} 2002; 132: 115--133.

\bibitem{XZ} 
Xia Yan, Zhang Guofen. Reliability equivalence functions in gamma distribution. 
{\it Appl. Math. Comput.} 2007; 187: 567--573. 
\end{thebibliography}
\end{document}